\documentclass[a4paper]{amsart}
\usepackage[T1]{fontenc}

\usepackage{longtable}
\usepackage{amssymb}

\usepackage[dvips, dvipsnames, usenames]{color}

\usepackage{color}
\usepackage{soul}

\usepackage{paralist}

\makeatletter

\newcommand{\fd}{finite-dimensional}

\newcommand\am{\mathbb A_m}

\newcommand\id{\operatorname{id}}

\newcommand{\ku}{\Bbbk}

\newcommand{\oc}{{\mathcal O}}
\newcommand{\Oc}{{\mathcal O}}

\newcommand\sn{\mathbb S_n}

\newcommand\bn{\mathbb B_n}

\providecommand{\\}{\\}

\numberwithin{equation}{section} 
\numberwithin{figure}{section} 
  \theoremstyle{plain}
  \newtheorem{thm}{Theorem}
  \theoremstyle{plain}
  
  \theoremstyle{remark}
  
  \theoremstyle{plain}
  
  \theoremstyle{plain}

\theoremstyle{remark}
\newtheorem*{acknowledgement*}{Acknowledgement}

\theoremstyle{definition}

\newcommand{\oct}{\mathfrak O}
\newcommand{\trid}{\triangleright}
\newcommand\toba{{\mathfrak B }}

\newcommand{\D}{{\mathcal D}}

\def\pf{\begin{proof}}
\def\epf{\end{proof}}

\begin{document}

\title[Nichols algebras over some sporadic simple groups]{Pointed Hopf algebras over\\ some sporadic simple groups}

\author[Andruskiewitsch, Fantino, Gra\~na, Vendramin]%
{N. Andruskiewitsch, F. Fantino, M. Gra\~na, L. Vendramin}

\thanks{Some of the results presented here are part of the PhD theses of FF and LV, work under the supervision of NA and MG, respectively.
\newline This work was partially supported by
 ANPCyT-Foncyt, CONICET, Ministerio de Ciencia y Tecnolog\'\i a (C\'ordoba) and Secyt (UNC)}

 \address{NA, FF: Facultad de Matem\'atica, Astronom\'\i a y F\'\i sica,
Universidad Nacional de C\'ordoba. CIEM -- CONICET. 
Medina Allende s/n (5000) Ciudad Universitaria, C\'ordoba,
Argentina}
\address{MG, LV: Departamento de Matem\'atica -- FCEyN,
Universidad de Buenos Aires, Pab. I -- Ciudad Universitaria (1428)
Buenos Aires -- Argentina}
\address{LV: Instituto de Ciencias, Universidad de Gral. Sarmiento, J.M. Gutierrez
1150, Los Polvorines (1653), Buenos Aires -- Argentina}

\email{(andrus|fantino)@famaf.unc.edu.ar}
\email{(matiasg|lvendramin)@dm.uba.ar}

\subjclass[2000]{16W30; 17B37}
\date{\today}
\begin{abstract}
Any finite-dimensional complex pointed Hopf algebra with group of group-likes
isomorphic to a sporadic group, with the possible exception of
the Fischer group $Fi_{22}$, the Baby Monster $B$ and the Monster $M$, is a group algebra.

\end{abstract}

\maketitle

\section{Introduction}
Let $\ku$ be an algebraically closed field of characteristic 0. In
this Note, we announce a new contribution to the classification of
finite-dimensional Hopf algebras over $\ku$.  As is known,
different classes of finite-dimensional Hopf algebras have to be
studied separately because the pertaining methods are radically
different. There is a method for pointed Hopf algebras (those
whose coradical is a group algebra $\ku G$) that has been applied
with satisfactory results when $G$ is abelian \cite{AS-05}; an
exposition of the method can be found in \cite{AS-cambr}.
Recently, it appeared that many finite simple (or almost simple)
groups $G$ admit very few \fd{}, pointed Hopf algebras with
coradical isomorphic to $\ku G$:

\begin{itemize}
    \item Any finite-dimensional complex pointed Hopf algebra with group of
        group-likes isomorphic to $\am$, $m\geq 5$, is a group
        algebra \cite{AFGV}.

    \item Same for the groups $SL(2,2^n)$, $n>1$ \cite{FGV}
        and $M_{20}$, $M_{21}=PSL(3,4)$ \cite{FGV2}.

     \item  Most of the
        pointed Hopf algebras over the symmetric groups have infinite
        dimension, with the exception of a short list of open possibilities,
        see \cite{AFGV, afz-2008} and references therein. More precisely, most of the irreducible Yetter-Drinfeld
        modules have infinite-dimensional Nichols algebras (see below).
\end{itemize}

This is a report on \fd{} pointed Hopf algebras over
sporadic simple groups. As part of our results, we have the
following.

\begin{thm}\label{teor:complete}
Let $G$ be any sporadic simple group, different from
the Fischer group $Fi_{22}$, the Baby Monster $B$ and the Monster $M$.
If $H$ is a
finite-dimensional pointed Hopf algebra with $G(H) \simeq G$, then
$H\simeq \ku G$.
\end{thm}

The Theorem holds more generally over any field of characteristic
0, since the property of being pointed is stable under extension
of scalars.

\subsection{Glossary}\label{subsec:glossary}
For the reader's convenience, we recall a few definitions that are central to our work. More information can be found
in \cite{AG1, AS-cambr}. Let $H$ be a Hopf algebra with comultiplication $\Delta$ and bijective antipode $\mathcal S$.

\smallbreak $\bullet$ An element $g\neq 0$ in $H$ is a \emph{grouplike} if $\Delta(g) = g\otimes g$; the set of all grouplikes is a group $G(H)$
with multiplication given by the product of $H$.

  \smallbreak $\bullet$ A \emph{Yetter-Drinfeld module} over $H$ is a left $H$-module $M$ that bears also a structure $\lambda: M \to H \otimes M$ of
  $H$-comodule, compatible with the action in an appropriate sense.
   If $H$ is finite-dimensional, then a Yetter-Drinfeld module is the same as a module over the Drinfeld double of $H$.
  For instance, if $H= \ku G$ is the group algebra of a finite group $G$,
  then a Yetter-Drinfeld module over $H$ is a left $G$-module $M$ that bears also a $G$-gradation $M = \oplus_{g\in G}M_g$,
  compatibility meaning that $h\cdot M_g = M_{hgh^{-1}}$ for all $h,g\in G$.

  \smallbreak $\bullet$ A \emph{rack} is a pair $(X,\trid)$ where $X$ is a non-empty set and $\trid:X\times X\to X$ is an operation such that
the map $\varphi_x = x\trid \underline{\quad}$
is bijective for any $x\in X$, and $x \trid(y\trid z) = (x\trid y) \trid(x\trid z)$
for all $x,y,z\in X$. A map $q:X\times
X\to GL(n,\ku)$ is a \emph{2-cocycle} of degree $n$ if $$q_{x,y\trid z}q_{y,z}=
q_{x\trid y,x\trid z}q_{x,z}, \quad\text{for all }x,y,z\in X. $$

  \smallbreak $\bullet$ A \emph{braided vector space} is a pair $(V,c)$ where $V$ is a vector space and $c\in GL(V \otimes V)$ fulfills the braid
  equation: $(c\otimes \id)(\id\otimes c)(c\otimes \id) = (\id\otimes c)(c\otimes \id)(\id\otimes c)$. Examples:

  (i) Any Yetter-Drinfeld module is a braided vector space in a natural way. 
  
  (ii). Let $X$ be a finite rack, $q$  a 2-cocycle of degree $n$, $V= \ku X\otimes\ku^{n}$, where $\ku X$ is the vector
space with basis $e_x$, $x\in X$. We denote $e_xv := e_x\otimes v$.
Consider the linear isomorphism $c^q:V\otimes V\to V\otimes V$,
$c^q(e_xv\otimes e_yw)=e_{x\trid y}q_{x,y}(w)\otimes e_xv$,
$x, y\in X$,  $v, w\in\ku^{n}$. Then $(V, c^q)$ is a braided vector space.
  
  The braided vector spaces
  arising as Yetter-Drinfeld modules over group algebras of finite groups can be presented in terms of racks and cocycles, see a
  bit more of information below.

  \smallbreak $\bullet$ We assume the reader familiar with the important notion of the
\emph{Nichols algebra} of a braided vector space, discussed at length in \cite{AS-cambr}.
In short, one of the possible definitions of the Nichols algebra $\toba(V)$ of a braided vector space
$(V,c)$ is as follows. Since $c$ satisfies the braid equation, it induces a representation of the braid group $\bn$,
 $\rho_n: \bn \to GL(V^{\otimes n})$, for each $n\ge 2$. Let $Q_n = \sum_{\sigma \in \sn} \rho_n(M(\sigma))\in End (V^{\otimes n})$,
 where $M: \sn \to \bn$ is the so-called Matsumoto section (not a morphism of groups, but preserves product when length is preserved).
Then $\toba(V)$ is the quotient of the tensor algebra $T(V)$ by $\oplus_{n\ge 2} \ker Q_n$, in fact a 2-sided ideal of $T(V)$.
If $c$ is the usual switch, then $\toba(V)$ is just the symmetric algebra of $V$; but in general the determination of a Nichols algebra
is quite a difficult task.

\section{Outline of the proof}
A complete proof of Theorem \ref{teor:complete} for the groups
$M_{22}$ and $M_{24}$ is contained in \cite{f}; the proof for the
other groups is included in \cite{AFGVe}.

We sketch now the proof in two main reductions.
The first one has been explained in several places, with detail in
\cite{AS-cambr}, but we include a brief summary for completeness.
We remind that if $U$ is a braided vector
subspace of $V$, then $\toba(U) \hookrightarrow \toba(V)$.

\subsection{A general reduction}\label{subsec:general}
Let $G$ be a finite group, $H$ a pointed Hopf algebra with
$G(H)\simeq G$. Then there are two basic invariants of $H$, a
Yetter-Drinfeld module $V$ over $\ku G$ (called the infinitesimal
braiding of $H$) and its Nichols algebra $\toba(V)$. We have
$\vert G\vert\dim\toba(V) \leq \dim H$. Therefore, the following
 statements are equivalent:

 \begin{asparaenum}[(1)]
     \item If $H$ is a \fd{} pointed Hopf
         algebra with $G(H) \simeq G$, then $H\simeq \ku G$.
     \item If $V\neq 0$ is a Yetter-Drinfeld module over $\ku G$,
         then $\dim \toba(V) = \infty$.
     \item If $V$ is an \emph{irreducible} Yetter-Drinfeld module over $\ku G$,
         then $\dim \toba(V) = \infty$.\label{it:viejo3}
 \end{asparaenum}

\subsection{Looking at subracks} We  focus  on (3) above.
The second reduction has been the basis of our recent papers. It
starts from the well-known classification of irreducible
Yetter-Drinfeld modules over $\ku G$ by pairs $(\Oc, \rho)$, where
$\Oc$ is a conjugacy class in $G$ and $\rho$ is an irreducible
representation of the stabilizer $G^s$ of a fixed point $s\in
\Oc$. Now, the definition of the Nichols algebra $\toba(\Oc,
\rho)$ of the corresponding Yetter-Drinfeld module $M(\Oc, \rho)$
just depends on the braiding. If $\dim \rho = 1$, then this
braiding depends only on the \emph{rack} $\Oc$ and a 2-cocycle $q:
\Oc \times \Oc \to \ku^{\times}$ \cite{AG1}. Namely, $\Oc$ is a
rack with the product $x\trid y := xyx^{-1}$, $M(\Oc, \rho)$ has a
natural basis $(e_x)_{x\in \Oc}$ and the braiding is given by
$c(e_x \otimes e_y) = q_{xy}e_{x\trid y}\otimes e_x$. If there
exists a subrack $X$ of $\Oc$ such that the Nichols algebra of the
braided vector space defined by $X$ and the restriction of $q$ is
infinite dimensional, then $\dim\toba(\Oc, \rho) = \infty$.

We recall some examples of racks which are relevant in this work.

\smallskip
\begin{asparaenum}[(i)]
    \item Abelian racks: those racks $X$ such that $x\trid y = y$ for all $x,y\in X$.

     \item\label{dos}  $\D_p$: the class of involutions
        in the dihedral group $\mathbb D_p$ (of order $2p$), $p$ a prime.

    \item\label{tres} $\oct$: the class of 4-cycles in  $\mathbb S_4$.

    \item Doubles of racks: if $X$ is a rack, then $X^{(2)}$ denotes the
        disjoint union of two copies of $X$ each acting on the other by left
        multiplication.
\end{asparaenum}

We are interested in finding subracks which are abelian, or
isomorphic to $\D_p^{(2)}$ or to $\oct^{(2)}$, by the following
reasons:

\smallskip
\begin{asparaenum}[(A)]
    \item If $X$ is abelian, then the corresponding braided vector space
        is of diagonal type. Braided vector spaces of diagonal type with
        \fd{} Nichols algebra where classified in \cite{H-all}; thus, we
        just need to check if the matrix $(q_{xy})$ belongs or not to the
        list in \cite{H-all}.
    \item\label{it:viejob}
        If $X$ is isomorphic either to $\D_p^{(2)}$ or to $\oct^{(2)}$, then
        for some specific cocycles, the related Nichols algebras have infinite
        dimension \cite[Ths. 4.7, 4.8]{AHS}.
\end{asparaenum}

\smallskip\emph{Variations.}
\begin{asparaenum}[(a)]
    \item If $\dim\rho >1$, similar arguments apply.
    \item Sometimes the rack $X$ is not abelian, but the braided vector space
        produced by $X$ and the $2$-cocycle can be realized with an abelian
        rack, by a suitable change of basis.
    \item Let $F< G$ be a subgroup, $s\in F$, $\Oc^F$, resp. $\Oc^G$ the
        conjugacy class of $s$ in $F$, resp. in $G$. If $\dim\toba(\oc^F,
        \tau) = \infty$ for any irreducible representation $\tau$ of $F^s$,
        then $\dim\toba(\oc^G, \rho) = \infty$ for any irreducible
        representation $\rho$ of $G^s$.
    \item A conjugacy class $\Oc$ is real if $\Oc = \Oc^{-1}$. It is quasireal
        if $\Oc=\Oc^m$ for some integer $m$, $1<m<N$, where $N$ is the order
        of the elements in $\Oc$. The search of subracks isomorphic to
        $\D_p^{(2)}$ or to $\oct^{(2)}$, as well as the verification that the
        restriction of the cocycle $q$ is as needed in \eqref{it:viejob}, is
        greatly simplified in a real (quasireal) conjugacy class \cite{AF}.
    \item We say that a rack $X$ is of type D if
	there exists  a decomposable subrack $Y = R\coprod S$ of $X$ such
	that $r\trid(s\trid(r\trid s)) \neq s$, for some $r\in R,\, s\in S$.
	If a conjugacy class $\oc$ is a rack of type D, then
	$\dim\toba({\oc,\rho})=\infty$ for any $\rho$
	(see \cite{AFGV} and Theorem 8.6 of \cite{HS}).

\end{asparaenum}

\subsection{Computations} We now fix a sporadic group $G$ as in Theorem
\ref{teor:complete}.  We extracted relevant information from the
\textsf{ATLAS} \cite{AtlasWeb} with the \texttt{AtlasRep} package \cite{AtlasRep}.
Then, we checked when a conjugacy class is real or
quasireal or of type D. We used \textsf{GAP} \cite{GAP} for the computations.

These tools allow to apply the techniques sketched above to all
pairs $(\Oc, \rho)$ and establish the validity of \eqref{it:viejo3}.

\subsection{} Some of these results were announced in several meetings:
\begin{itemize}
	\item Hopf Algebras and Related Topics, A conference in honor of
		Professor Susan Montgomery. University
		of Southern California, Los Angeles, USA. February 2009.
	
	\item IV Encuentro Nacional de \'Algebra, Córdoba, Argentina. August,
		2008.
	
	\item First De Brún Workshop on Computational
		Algebra, National University of Ireland, Galway, Ireland.
		August, 2008
	
	\item Groupes quantiques
		dynamiques et cagories de fusion. CIRM, Luminy, France. April
		2008.
\end{itemize}


\newcommand{\etalchar}[1]{$^{#1}$}

\end{document}